
\documentstyle{amsppt}

\topmatter
\title Dynkin Graphs and Triangle Singularities
\endtitle
\author
Tohsuke Urabe
\endauthor\affil
Department of Mathematics\\
Tokyo Metropolitan University\\
Minami-Ohsawa 1-1,\ Hachioji-shi\\
Tokyo 192-03 \ Japan\\
(E-mail: urabe\@math.metro-u.ac.jp)\endaffil
\subjclass 14B12, 14J17, 32S30\endsubjclass
\endtopmatter

\document
\hyphenation{co-ef-fi-cient}
In Arnold's classification list of singularities (Arnold \cite{1}.)
we find interesting singularities to be studied. Though we find
singularities of any dimension in Arnold's list, we consider
singularities of dimension two in particular.  Among them there is a
class called exceptional singularities or triangle singularities.
This class consists of fourteen singularities.  It is known that
they are closely related to K3 surfaces with the structure of
elliptic surfaces.  (Looijenga \cite{4}.)  Here we would like to
consider the following nine singularities of these fourteen ones:

$$\matrix
E_{\,12},&Z_{\,11},&Q_{\,10},\\
  &&\\
  E_{\,13},&Z_{\,12},&Q_{\,11},\\
  &&\\
  E_{\,14},&Z_{\,13},&Q_{\,12}.\endmatrix$$
(The remaining five triangle singularities are
$W_{\,12},\;W_{\,13},\;S_{\,11},\;S_{\,12}$
and
$U_{\,12}$.)  We assume that the ground field is the complex field
${\bold C}$.

Recall here that a connected Dynkin graph of type $A$, $D$  or
$E$ corresponds to a surface singularity called a rational double
point.
Let  $\Xi$   be a class of surface singularities.  By  $PC(\Xi)$ we
denote the set of Dynkin graphs $\Gamma$ with several components
such that there exists a small deformation fiber  $Y$  of a
singularity belonging to  $\Xi$  satisfying the following
conditions:

\roster
\item  $Y$  has only rational double points as singularities.
\item The combination of rational double points on  $Y$  corresponds
exactly to  $\Gamma$.  (The type of each component of  $\Gamma$
corresponds to the type of the singularity on  $Y$  and the number
of components of each type corresponds to the number of
singularities of each type on  $Y$.)
\endroster

Note that by definition every graph in $PC(\Xi)$ has only components
of type $A$, $D$ or $E$.  We would like to study $PC(\Xi)$ for
$\Xi=E_{\,12},\;Z_{\,11},\;...,\;Q_{\,12}$

\midspace{5mm}\pagebreak
\proclaim{Theorem}
Let  $\Xi$  be one of the
above  nine classes of singularities.   The following two conditions
are equivalent.

\roster
\item"(A)"  $\Gamma\in PC (\Xi)$.

\item"(B)"  The Dynkin graph  $\Gamma$  has only  components of type
$A$,
$D$ or $E$, and can be made from the essential basic graph depending
on  $\Xi$ by a combination of two of elementary transformations and
tie transformations. \endroster

$$\vbox{
\hbox{\qquad
The respective essential basic graph
}
\hbox{\it
corresponding to the above nine singularities
}
}$$
$$\matrix
E_{\,8},&E_{\,7},&E_{\,6}\\
  &&\\
  E_{\,8}+BC_{\,1},&E_{\,7}+BC_{\,1},&E_{\,6}+BC_{\,1}\\
  &&\\
  E_{\,8}+G_{\,2},&E_{\,7}+G_{\,2},&E_{\,6}+G_{\,2}
\endmatrix$$
\endproclaim
\midspace{0pt}

\rm
In the above list of essential basic graphs the plus sign  $+$
denotes the disjoint union of graphs.
In the above condition (B) an elementary transformation and
a tie transformation are operations by which we can make a new
Dynkin graph from a given Dynkin graph.  We give the definition of
them below.  In the condition (B) four kinds of combinations --
``elementary" twice, ``tie" twice, ``elementary" after ``tie", and
``tie" after ``elementary" -- are all permitted.

\proclaim{ Definition}\rm
 (An elementary transformation)\quad The following
procedure  is called an elementary transformation of a Dynkin graph.

\roster
\item Replace each connected component by the corresponding extended
Dynkin graph.

\item Choose in arbitrary manner at least one vertex from each
component (of the extended Dynkin graph) and then remove these
vertices together with edges issuing from them.
\endroster\endproclaim

We can find the definition of the extended Dynkin graph in any book
on Lie algebras. (Bourbaki \cite{2}.) They can be made by adding one
vertex and one or two edges to each connected component of the
Dynkin graph. The position of the added vertex and edges depends on
the type of the component.  Also we can find the definition of the
coefficients of the maximal root in any book on Lie algebras.

\proclaim{ Definition}\rm
(A tie transformation)\quad Assume that applying the following
procedure to a Dynkin graph  $\Gamma$ , we have obtained the Dynkin
graph
${\overline\Gamma}$. Then we call the following procedure a tie
transformation of a Dynkin graph.

\roster
\item Add one vertex and a few edges to each component
of  $\Gamma$  and make it into the extended Dynkin graph of the
corresponding type. Moreover attach the corresponding coefficient of
the maximal root to each vertex.

\item Choose in an arbitrary manner subsets  $A$, $B$  of the set of
the vertices of the extended graph $\tilde \Gamma $ satisfying the
following conditions:\endroster

\hbox{\qquad\quad\hfill\vbox{\hsize=11.5cm \roster
\item"$\langle a\rangle$" $A\cap B=\emptyset $.

\item"$\langle b\rangle$" Let $V$  be the set of vertices of an
arbitrarily chosen component  $\tilde \Gamma '$
of  $\tilde \Gamma $.  Let  $\ell $
be the number of elements in  $V\cap A$ and
$n_{\,1},\;n_{\,2},\;...,\;n_{\,\ell }$ be the numbers attached
to $V\cap A$.  Furthermore let  N  be the sum of numbers attached to
$V\cap B$.  (If  $V\cap B=\emptyset $, then  $N=0$ .)  Then the
greatest common divisor of the  $\ell +1$ numbers
$N,\;n_{\,1},\;n_{\,2},\;...,\;n_{\,\ell }$ is 1.
\endroster}}
\roster
\item[3] Erase all attached integers and remove
vertices belonging to  $A$  together with edges issuing from them.

\item[4] Draw another new vertex  $\bigcirc$  corresponding to a
root
$\alpha$  with $\alpha ^2=2$. Connect this new vertex  $\bigcirc$
and each vertex in  $B$  by an edge.  \endroster\endproclaim

\proclaim{Remark}\rm
Often the resulting graph  $\bar \Gamma $ after the above
procedure (1) -- (4) is not a Dynkin graph. We consider only the
cases where the resulting graph $\bar \Gamma $ is a Dynkin graph and
then we call the above procedure a tie transformation.  Under this
restriction the number $\#(B)$ of elements in the set  $B$  satisfies
$0\le \#(B)\le 3$.  $\ell =\#(V\cap A)\ge 1$.
\endproclaim

Here we give some explanation on Dynkin graphs and root systems of
type $BC$.  A root system  $R$  is a finite subset of a Euclidean
space satisfying axioms on symmetry.  Usually we assume moreover the
following axiom (*) of the reduced condition:

\qquad\qquad\qquad\qquad\qquad\qquad
If  $\alpha \in R$, then
$2\alpha \notin
R$\qquad\qquad\hfill\hfill(*)

\noindent
Under these axioms we obtain irreducible root systems of type $A$,
$B$, $C$, $D$, $E$, $F$  and  $G$  as in any book on Lie algebras.
However, under the absence of the axiom (*) we have further a series
of irreducible root systems, which are called of type $BC_k$
$\left( {k=1,\;2,\;3,\;...} \right)$. (Bourbaki \cite{2}.)  It is
easy to generalize the concept of Dynkin graphs to root systems of
type
$BC$.  (Urabe \cite{6}.)  The Dynkin graph of type $BC_1$ is the
following: $\bigotimes$

We explain the meaning of this $BC_1$ graph.  Recall first the
meaning of Dynkin graphs.  Let  $R$  be an irreducible root system
and $\Delta \subset R$ be the root basis.  We can assume that the
longest root $\alpha\in R$ satisfies $\alpha ^2=2$ after normalizing
the inner product of the ambient Euclidean space.  The Dynkin graph
$\Gamma$ of  $R$  is the graph drawn by the following rules:  (1)
The vertices  of $\Gamma$ have one-to-one correspondence with the set
$\Delta$ (the root basis).  (2) Two vertices in $\Gamma$
corresponding to two elements $\alpha ,\;\beta \in \Delta $ are
connected by an edge in $\Gamma$ if and only if the inner product
$\left( {\alpha ,\;\beta } \right)\ne 0$.

If  $R$  is of type  $A$, $D$  or  $E$,  then  $R$  consists of only
roots $\alpha$ with $\alpha ^2=2$,  and every $\alpha \in \Delta $
satisfies $\alpha ^2=2$.  Therefore in these cases every vertex in
the Dynkin graph can be denoted by a small white circle $\bigcirc$.

If  $R$  is of type $BC_1$, then $\Delta$ consists of a unique root
$\delta$ with $\delta ^2=1/ 2$ and $R=\left\{ {\,-2\delta
,\;-\delta ,\;\delta ,\;2\delta \,} \right\}$. The vertex in the
Dynkin graph corresponding  to a root $\delta$ with $\delta ^2=1/
2$ is denoted by $\bigotimes$. The $BC_1$ graph is the graph
consisting of a unique vertex of this kind.  In this case the
maximal root $\eta $ is equal to $2\delta$, and thus the extended
Dynkin graph, i.e., the graph corresponding to $\Delta ^+=\Delta
\cup \left\{ {-\eta } \right\}$  is the
following:\phantom{i}\qquad\qquad\qquad\qquad\qquad\qquad (The edge
is bold. The numbers are the coefficients of the maximal root.)

If  $R$  is of type $G_2$, then  $\delta$ consists of two elements
$\alpha$  with $\alpha^2=2$ and $\gamma$ with $\gamma^2=2/3$.
We denote  the  vertex  corresponding  to a  root
$\gamma$   with  $\gamma^2=2/3$ by
\phantom{i}\qquad.
Our Dynkin graph of type $G_2$ is the following;
\qquad\qquad\qquad\qquad\quad
and our extended Dynkin graph of
type $G_2$ is the following (The numbers are the coefficients of the
maximal root.):

\midspace{0pt}

Note that as a result of an elementary or a tie transformation, a
graph  consisting of a unique vertex corresponding to $\gamma$ with
$\gamma^2=2/3$ can appear.  We call the graph
\qquad\quad the Dynkin graph of type $G_1$. This
corresponds to the root system
$R=\left\{ {\,-\gamma ,\;\gamma \,} \right\}$ with $\gamma^2=2/3$.
The extended Dynkin graph of type $G_1$ is the following:
\phantom{i}\qquad\qquad\qquad\qquad\qquad\qquad.
(The edge is bold.
The numbers are the coefficients of the maximal root.)

We can explain why we do not use the standard
expression \phantom{i}\qquad\qquad\qquad\qquad
of the $G_2$ graph. (Bourbaki \cite{2}.)  If we
use the standard expression, we cannot define the concept of the
$G_1$ graph.

Now, for an irreducible root system  $R$ of the remaining types,
i.e., of type $B_k$, $C_k$, $F_4$ or $BC_l$ with $l\geq 2$, the root
basis $\Delta$ contains a root $\beta \in \Delta $ with
$\beta^2=1$. However, in the case of our nine triangle singularities
such a root $\beta$ with $\beta^2=1$ never appears.  Therefore in
our case the type of a connected Dynkin graph is either $A_k$ with
$k\geq 1$, $D_l$ with $l\geq 4$, $E_6$, $E_7$, $E_8$, $G_2$, $G_1$
or $BC_1$.

Note that since we have assumed that the Dynkin graph $\Gamma$ in our
Theorem has only components of type  $A$, $D$  or  $E$ ,  any Dynkin
graph with a component of type  $G_2$, $G_1$ or $BC_1$ made by two
transformations has no meaning, and is to be thrown out.

\proclaim{Example}\rm
We show $A_{\,7}+A_{\,4}\in PC (Z_{\,13})$ and
$D_{\,8}+A_{\,2}\in PC (Z_{\,13})$.

For $Z_{13}$ the corresponding essential basic graph is
$E_{\,7}+G_{\,2}$. We can start from $E_{\,7}+G_{\,2}$.
As an example, we apply a tie transformation to this graph.  After
the first  step of the transformation the following graph is
obtained and we can choose the subsets  $A$  and  $B$  as follows.

\hbox to \hsize{}
\midspace{40mm}

Obviously the condition $\langle a\rangle$ $A\cap B=\emptyset $ is
satisfied. For the component $E_7$, $\ell =\#(V\cap A)=1$ and
$n_1=1$, $N=1$. Thus $G .\,C .\,D .(n_{\,1},\;N)=1$. For the
component $G_2$, $\ell =1$, $n_1=1$, $N=0$ and $G .\,C .\,D
.(n_{\,1},\;N)=1$. The condition $\langle b \rangle$ is also
satisfied.  In the next step all vertices in  $A$  are erased, and
drawing a new vertex $\bigcirc$, we connect it and the vertex in
$B$  by an edge.  One knows that the resulting graph is
$E_{\,8}+G_{\,2}$.

We can apply a transformation once more starting from
$E_{\,8}+G_{\,2}$.  First we apply a tie transformation.

\midspace{50mm}

\noindent
The above choice of $A$  and  $B$  satisfies the conditions and
we get the graph  $A_7+A_4$  as the result.  By our Theorem
$A_{\,7}+A_{\,4}\in PC (Z_{\,13})$.

If we apply an elementary transformation to  $E_8+G_2$, and if
we erase two vertices as follows, we obtain the graph  $D_8+A_2$.

\hbox to \hsize{}
\midspace{50mm}

By our Theorem one can conclude  $D_{\,8}+A_{\,2}\in PC (Z_{\,13})$.
\endproclaim

Below we sketch the verification of our Theorem briefly.

First we apply the results in Looijenga \cite{4}.  In \cite{4}
Looijenga shows  that our singularity is closely related K3 surfaces
and that by the theory of periods for K3 surfaces we can reduce our
problem into a problem on the lattice theory.  Let  $\Lambda_N$  be
the even unimodular lattice with signature  $\left( {16+N,\;N}
\right)$. It is unique up to isomorphisms  if  $N\geq 1$.  A certain
lattice
$P$   is defined corresponding to each $\Xi$ of the nine triangle
singularities.  We can consider the quotient module
$\Lambda_N/P$  when a primitive embedding  $P\subset \Lambda
_{\,N}$ is given.  A bilinear form on  $\Lambda_N/P$  with values
in rational numbers can be defined.  By Looijenga one knows that
for a Dynkin graph  $\Gamma$ with only components of type  $A$, $D$
or $E$, $\Gamma \in PC (X)$ if and only if the associated root
lattice  $Q(\Gamma )$ has a lattice embedding into  $\Lambda
_{\,3}/ P$ satisfying certain conditions.  (Note that the suffix
of  $\Lambda$  is $N=3$  here.)

Second we translate Looijenga's conditions on the lattice theory
into a simpler condition.  We say that an embedding $Q(\Gamma
)\subset \Lambda _{\,N}/ P$ is {\sl full} if the root system of
$Q(\Gamma )$  and the root system of the {\sl primitive hull} of
$Q(\Gamma )$ in  $\Lambda_N/P$ coincide.  (For a submodule  $M$  of
a module  $L$, the set  $\tilde M=\{\,x\in L\,\mid\, \hbox{\rm For
some non-zero integer } m,$ $ mx\in M\}$ is called the
primitive hull of $M$  in  $L$.  Obviously  it is a submodule
containing
$M$.)  We consider here root systems including roots  $\delta$
with  $\delta ^{\,2}=1/ 2$ and roots $\gamma$ with  $\gamma^2=2/3$.
One knows that a lattice embedding  $Q(\Gamma )\subset \Lambda
_{\,3}/ P$ satisfies Looijenga's conditions if and only if it is
full.

Let  $\overline {PC }(\Xi)$ denote the set of Dynkin graphs
satisfying the condition (B) in our Theorem.

The inclusion relation $\overline {PC }(\Xi)\subset PC (\Xi)$
is an immediate consequence of our general theory of elementary
transformation and tie transformations.  Indeed, let  $\Gamma'$  be
a Dynkin graph obtained from a Dynkin graph  $\Gamma$ by an
elementary transformation or a tie transformation. We can show that
if there exists a full embedding $Q(\Gamma )\subset \Lambda _{\,N}/
P$,  then there exists a full embedding $Q(\Gamma ')\subset \Lambda
_{\,N+1}/ P$.  Note here that the suffix of $\Lambda$ increases by
one.  Besides, for some primitive embedding $P\subset \Lambda
_{\,1}$ the corresponding essential basic graph $\Gamma_0$ has a
full embedding $Q(\Gamma _0)\subset \Lambda _{\,1}/ P$.  Thus we
can conclude $\overline {PC }(\Xi)\subset PC (\Xi)$.

Let $\mu$ be the Milnor number of one of the nine triangle
singularities under consideration.  This number is equal to the
suffix of the corresponding symbol of the singularity.  (For
$E_{12}$,  $Z_{12}$  and  $Q_{12}$  $\mu=12$.)  Let $\Gamma$ be a
Dynkin graph with  $r$  vertices.  We can show that if $\Gamma \in
PC (\Xi)$,  then $r\leq \mu -2$.  Besides, if  $r\leq \mu -5$,
then conditions $\Gamma \in PC (\Xi)$ and $\Gamma \in \overline
{PC }(\Xi)$ are equivalent.  The last assertion follows from
Meyer's theorem ``Any indefinite rational quadratic form represents
zero, if the number of variables is greater than or equal to five."

In order to show the opposite inclusion relation $\overline {PC
}(X)\supset PC (X)$ it suffices to show $\left( {S_{\,\mu }\cap
\overline {PC }(\Xi )} \right)\cup \left( {S_{\,\mu }-PC (\Xi )}
\right)=S_{\,\mu }$,  where $S_{\,\mu }$ denotes the set of Dynkin
graphs $\Gamma$ with only components of type $A$, $D$  or  $E$
whose number  $r$  of vertices satisfies $\mu -4\le r\le \mu -2$.
To tell the truth, we could not succeed in finding any effective
method to show this equality except case-by-case checking.  This is
a weak point of our theory.  I regret this fact and hope that
somebody can improve it.  The theory of monodromy groups of elliptic
surfaces may be effective for the improvement.  Anyway, by the
elementary lattice theory,  the surface theory in the algebraic
geometry, and the $p$-adic lattice theory due to Nikulin (Nikulin
\cite{5}), we can accomplish the checking.

Details of the verification will appear elsewhere.

Before concluding this article we would like to refer to the
remaining five hypersurface triangle singularities $W_{12}$,
$W_{13}$, $S_{11}$, $S_{12}$ and $U_{12}$.  Recall here in
particular that the number of transformations in the condition (B)
in our Theorem for the nine singularities is two. For the remaining
five singularities $W_{12}$, $W_{13}$, $S_{11}$, $S_{12}$ and
$U_{12}$, if we try to formulate the corresponding theorem including
a description in which the number of transformations is two, the
formulation becomes very complicated and it is not worth
mentioning.  This is because the property of the {\sl Milnor
lattice} for $W_{12}$, $W_{13}$, $S_{11}$, $S_{12}$ and $U_{12}$ is
very different from that of the nine singularities in this article.

 We consider one of fourteen hypersurface triangle singularities and
let   $F$  be the corresponding Milnor fiber.  The pair
$(L,\;-(\;\;,\;))$ of the second homology group  $L=H_{\,2}(F,\,{\bold
Z})$ of  $F$  and $(-1)$ times the intersection form $(\;\;,\;)$ is
called the Milnor lattice.  Let $H={\bold Z}u+{\bold Z}v$ denote the
hyperbolic plane, i.e., a lattice of rank $2$ with $u^2=v^2=0$ and
$\left( {u,\,v} \right)=1$.

 For any of fourteen cases  $L$  has the following decomposition

\qquad\qquad\qquad\qquad\qquad\qquad
$L\cong M\oplus H\oplus H$\qquad\qquad\hfill\hfill(**)

\noindent
where  $M$  is a positive definite lattice, and  the symbol
$\oplus$  denotes the orthogonal direct sum.  For the nine
singularities considered in this article for every decomposition
(**) the co-root system $R^\vee $ of  $M$,  i.e., the set
$$R^\vee =\{\,x\in M\,\mid\,x^2=2,\;4 \hbox{\rm \ or\ } 6. \hbox{\rm
\ For every\ } y\in M\   2(x,\,y)/x^2 \hbox{\rm \ is an integer.}\}$$
spans  $M$  over  ${\bold Q}$.  However, for $W_{12}$, $W_{13}$,
$S_{11}$, $S_{12}$ and $U_{12}$ for every decomposition (**) the
root system never spans  $M$.

Instead of a theorem with the number of transformations two,
we can formulate a theorem in which the number of transformations is
one.  In this case the basic graph used at the start of
transformations is not Dynkin but a so-called Gabri\'{e}lov graph.
(Gabri\'{e}lov \cite{3}.)  For $W_{12}$, $W_{13}$, $S_{11}$, $S_{12}$
and
$U_{12}$ under this formulation we can get theorems worth
mentioning.  Besides, even for the nine singularities considered in
this article also theorems under the formulation with the number of
transformations one are worth mentioning.

These results will appear elsewhere.

\enddocument